# Multi-Agent Receding Horizon Games Framework for Autonomous Market Participation

**Parth Brahmbhatt[a] and Styliani Avraamidou[a*]**

[a] University of Wisconsin-Madison, Department of Chemical and Biological Engineering, Madison, WI, USA
* Corresponding Author: avraamidou@wisc.edu

## ABSTRACT

Electrochemical distributed energy resources (DERs), such as electrolyzers and battery energy storage systems, consume large amounts of electricity. Buying that power directly from wholesale markets can sharply reduce their operating costs. Accessing these markets, however, requires a minimum bid size that individual small or medium-scale units struggle to meet efficiently. Grouping several units so they bid as a single participant clears this barrier. It also lets each unit keep using its own low-cost renewable generation. The common remedy is to hire a third-party aggregator that bids on behalf of the group. However, aggregators charge commissions. They also solve a single centralized optimization whose outcome often favors some units while leaving others with little benefit. A fairer alternative is peer-to-peer (P2P) participation. Here, units form a self-governing group and bid jointly with no central authority. Each unit can then rely on its own renewables, while the group collectively meets the market minimum threshold. In practice, each unit is an independent, self-interested entity that is unwilling to reveal its data. This is why such coordination is naturally posed as a game-theoretic distributed optimization problem, rather than a single joint optimization. Existing P2P methods, however, address only single-round spot markets. They ignore the two-stage structure of real wholesale markets, where participants commit a day ahead and then continuously adjust in real time. We close this gap with a two-stage receding-horizon generalized Nash equilibrium (GNE) game. Each unit re-optimizes its strategy every five minutes over a rolling one-hour horizon, while the group satisfies both market stages collectively. The units exchange only publicly visible aggregate power, never private cost or production data. We apply the approach to a six-unit fleet on the PJM market. It delivers 37% higher profit, 29% lower electricity cost, and 92% less renewable curtailment than individual participation, with every unit better off.

## INTRODUCTION

Electrochemical processes such as water electrolyzers and battery energy storage systems (BESS) are a growing class of distributed energy resources (DERs) that consume large amounts of electricity. Electricity cost accounts for 60-80% of their total operating cost [1]. Wholesale electricity markets offer significantly cheaper power during periods of low demand or excess renewable supply, making direct market access a powerful tool for reducing operating costs. However, wholesale markets require participants to commit a minimum power level. In PJM, this threshold is 100 kW [2]: a unit either stays fully out of the market or must purchase at least 100 kW from the grid during any committed interval. This same floor applies in real-time; any upward adjustment above the day-ahead plan must also meet the 100 kW minimum. Small units cannot meet this threshold alone and are excluded from the market entirely. Medium-scale units (125-250 kW) can meet it individually, but still face a hidden problem. Markets operate in two rounds: units first submit a plan the day before (day-ahead), then make real-time adjustments every five minutes. A rule governing these real-time adjustments requires that any increase from the day-ahead plan must also meet the 100 kW minimum. For a unit with co-located solar or wind generation, this creates two distinct problems. When renewable output is sufficient to power the unit entirely, the solution is simple: the unit exits the market and consumes its own generation. But when renewable output is partial enough to reduce grid needs below 100 kW but not enough to run the unit alone, the market floor forces the unit to purchase at least 100 kW from the grid regardless, displacing renewable generation that must then be curtailed. Medium-scale units are therefore individually capable of market participation but individually worse off for doing so during renewable generation peaks.

The typical solution is to hire a third-party aggregator who bids on behalf of a group of units [3]. Aggregators solve complex centralized optimization problems, charge commissions, and often make decisions that strongly favor some units while providing little benefit to others [4]. A more equitable alternative is peer-to-peer (P2P) market participation, where units form a self-governing group and bid jointly without any central authority [5-7]. In practice, these units are independent, self-interested entities, often under separate ownership, and are unwilling to disclose their private cost structures or production data to one another. Their coordination, therefore, cannot be cast as a single centralized optimization; it is instead naturally posed as a game among the units, which is why P2P participation is typically formulated through game-theoretic distributed optimization. Each unit in such a group optimizes its own profit while the group as a whole satisfies market rules, so renewable-rich units can step back from the market when generation is high, while grid-only units cover the group's minimum commitment. Existing P2P frameworks [5-7, 15, 16], however, only address single-round spot markets and do not account for the sequential structure of real wholesale markets, where units must first commit to a plan a day ahead and then continuously adapt that plan as prices and generation change in real time.

This paper fills that gap. We propose a receding-horizon two-stage game in which each unit independently re-optimizes its strategy every five minutes over a rolling one-hour horizon, while the group collectively satisfies both rounds of market rules. Each unit solves only its own local problem, no unit shares its cost structure or production data with others, the only information exchanged between units is the aggregate power consumed by the group, which is already publicly visible to all market participants, and a fully distributed coordination algorithm drives the group toward a Generalized Nash Equilibrium (GNE) [8-10]: a state where no unit can improve its own profit by changing its strategy unilaterally. Combining a two-stage wholesale market structure with receding-horizon replanning within such a peer game has not been previously reported for electrochemical process units.

# PROBLEM FORMULATION

## Two-Stage Market Participation

We model a wholesale market that clears in two stages. Time is discretized into steps of length $t$ [h]. In the *day-ahead* (DA) stage, run once per day, each participant submits a schedule of grid-import quantities and an hourly on/off commitment over the $T$ steps of the coming day (indexed by $t$). In the *real-time* (RT) stage, the schedule is adjusted to track actual prices and renewable generation; we solve it in receding-horizon fashion, re-optimizing over a short $H$ look-ahead window (indexed by $k = 0, \dots, H-1$) at each step but committing only the first step before advancing [11]. Figure 1 illustrates this two-stage structure: the DA stage runs once at the start of the day to fix the DA commitment schedule, while the RT stage re-optimizes repeatedly in a rolling fashion, advancing one step at a time within each hour-long look-ahead window.

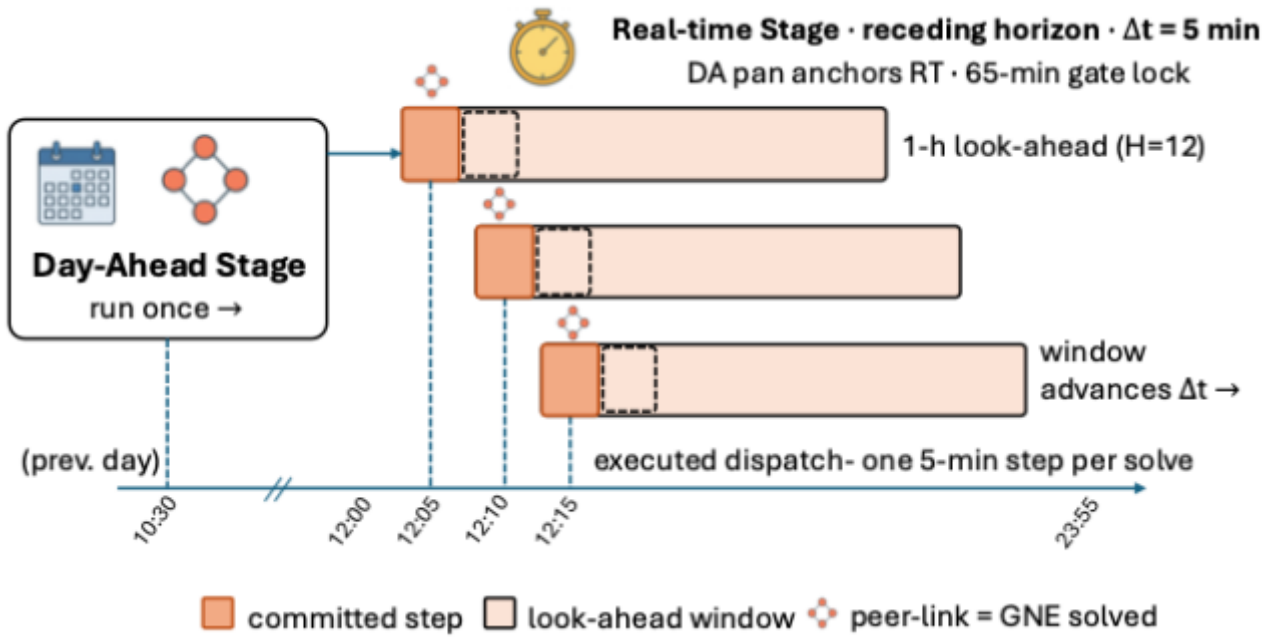

**Figure 1.** Two-stage market operation timeline.

Both stages obey a minimum-bid rule: a participant either stays out ($p = 0$) or imports at least $L_{min}$ [kW]. With a binary commitment $u \in \{0,1\}$ and rated power $P^{max}$ [kW],

$$L_{\min}\, u \ \le\ p \ \le\ P^{\max}\, u. \tag{1}$$

A second rule links the stages: in RT, any increase of import above the DA quantity must itself be at least $L_{min}$. We state both rules first for a single asset, then for a coalition.

## Asset Model

Let $\mathcal{N} = \{1, \dots, N\}$ index the fleet. Each asset $i$ comprises of a continuously-modulating electrochemical device (rated $P_i^{max}$, yield $\eta_i$), and/or optionally co-located renewable generation and a battery. The decision variables at step $t$ are: grid import $p_{i,t} \ge 0$, converter load $e_{i,t}$, renewable curtailment $c_{i,t}$, battery charge/discharge $p_{i,t}^c, p_{i,t}^d$ with mode binary $\delta_{i,t}$ and state of charge $E_{i,t}$, and the commitment $u_{i,t}$. A power balance routes renewable generation $g_{i,t}$, battery flow and grid import into the converter,

$$e_{i,t} = p_{i,t} + g_{i,t} - c_{i,t} + p_{i,t}^d - p_{i,t}^c, \tag{2}$$

subject to the operating limits

$$0 \le e_{i,t} \le P_i^{\max}, \qquad 0 \le c_{i,t} \le g_{i,t}, \tag{3}$$

$$|e_{i,t} - e_{i,t-1}| \le R_i\, \Delta t, \tag{4}$$

$$p_{i,t}^c \le \overline{P}_i^c\, \delta_{i,t}, \quad p_{i,t}^d \le \overline{P}_i^d\, (1 - \delta_{i,t}), \tag{5}$$

$$E_{i,t+1} = E_{i,t} + \eta^c p_{i,t}^c \Delta t - p_{i,t}^d \Delta t / \eta^d, \quad 0 \le E_{i,t} \le \overline{E}_i. \tag{6}$$

Here $R_i$ is the ramp limit, $\overline{P}_i^c, \overline{P}_i^d$ the battery power limits, $E_i$ its capacity, and $\eta^c, \eta^d$ its efficiencies; (5) forbids simultaneous charge and discharge. Asset $i$ produces $\eta_i e_{i,t} \Delta t$ units of product per step. Let $\mathcal{L}_i$ denote the local feasible set (2)-(6).

## Day-Ahead Planning Problem

In the DA stage, asset $i$ minimizes net cost, electricity purchased at the forecast price $\lambda_t$ [$/MWh] minus product revenue at rate $r_i$ [$/unit] over the full day:

$$\min_{(p_i, e_i, u_i, \dots) \in \mathcal{L}_i} \sum_{t=1}^{T} \left( \frac{\lambda_t}{1000}\, p_{i,t} - r_i\, \eta_i\, e_{i,t} \right) \Delta t \tag{7}$$

subject to the asset model (2)-(6), the semi-continuous bid (1) at every step, and a daily product-demand floor with target $D_i$

[units]:

$$\sum_{t=1}^{T} \eta_i\ e_{i,t}\,\Delta t\ \ \geq\ \ D_i. \tag{8}$$

The solution yields the DA commitment $u_i^{\mathrm{DA}}$ and import schedule $p_i^{\mathrm{DA}}$, which anchor the RT stage.

## Real-Time Receding-Horizon Problem

At each RT step, the asset re-optimizes over the next $H$ steps, subject to the same asset model (2)-(6). The objective adds a quadratic penalty, weighted by $\gamma_i \geq 0$, that discourages costly departures from the committed DA position $p_{i,k}^{\mathrm{DA}}$:

$$\min_{(p_i,e_i,\dots)\in\mathcal{L}_i} \sum_{k=0}^{H-1} [\Big(\frac{\lambda_k}{1000} p_{i,k} - r_i \eta_i e_{i,k}\Big)\Delta t + \frac{\gamma_i}{2}\big(p_{i,k} - p_{i,k}^{\mathrm{DA}}\big)^2] \tag{9}$$

The commitment $u_i$ is locked at its DA value $u_i^{\mathrm{DA}}$ (so the minimum-bid rule (1) holds with $u$ fixed), and the demand floor (8) becomes a cumulative catch-up target over the remaining steps with a heavily penalized slack, so a single poor forecast never renders the step infeasible. Only step $k = 0$ is executed; the window advances.

## Coalition Coupling and the Generalized Nash Equilibrium

When the assets form a coalition, they connect to the grid through a single metered point and bid as one market participant, as illustrated in Figure 2. The minimum-bid rule (1) and the two-stage increment rule then apply to the aggregate import, not to any individual asset ($L_{\max}$ is the aggregate import cap, $\infty$ if uncapped):

$$\sum_{i\in\mathcal{N}} p_{i,t}\ \in\ \{0\} \cup [L_{\min},\, L_{\max}], \tag{10}$$

$$\sum_{i\in\mathcal{N}} \big(p_{i,k} - p_{i,k}^{\mathrm{DA}}\big)\ \in\ \{0\} \cup [L_{\min},\, \infty). \tag{11}$$

This aggregate form is the source of the coalition's advantage: a renewable-rich asset may set $p_{i,t} = 0$ and consume its own generation while grid-only partners jointly supply the floor $L_{\min}$. Each asset, however, still minimizes only its own cost $J_i$ ((7) or (9)). Because the coupling constraints (10-11) make each asset's feasible set depend on the others' decisions $\mathbf{p}_{-i}$, the appropriate solution concept is a GNE [8, 9]: a profile $\mathbf{p}^\star = (p_1^\star, \dots, p_N^\star)$ at which no asset can lower its own cost by changing its strategy alone,

$$p_i^\star \in \arg\min_{p_i\,:\,(p_i,\mathbf{p}_{-i}^\star)\in\mathcal{C}} J_i(p_i), \qquad \forall\, i \in \mathcal{N}, \tag{12}$$

where $\mathcal{C}$ collects the shared constraints (10-11). Since the asset costs are separable and couple only through the shared import limits, the game is an exact potential game and its variational GNE coincides with the coalition optimum [5].

## Distributed Solution via MIQP-ADMM

We compute the GNE with a distributed alternating direction method of multipliers (ADMM) in which each asset solves only its own mixed-integer quadratic program (MIQP) and the only quantity exchanged is the aggregate import [12, 13]. Let $z_i$ be a consensus copy of $p_i$, $\mu_i$ the associated dual price, and $\rho > 0$ the ADMM penalty parameter. Each iteration $\ell$ performs:

$$p_i^{\ell+1} = \arg\min_{p_i\in\mathcal{L}_i} J_i(p_i) + \mu_i^{\ell\top} p_i + \frac{\rho}{2}\big\|p_i - z_i^{\ell}\big\|_2^2, \tag{13}$$

$$z^{\ell+1} = \Pi_{\mathcal{C}}\Big(\sum_i \big(p_i^{\ell+1} + \mu_i^{\ell}/\rho\big)\Big), \tag{14}$$

$$\mu_i^{\ell+1} = \mu_i^{\ell} + \rho\,\big(p_i^{\ell+1} - z_i^{\ell+1}\big). \tag{15}$$

The local MIQP (13) retains all of asset $i$'s integrality (commitment $u_i$, battery mode $\delta_i$) and is solved independently in parallel. The consensus step (14) couples the assets by projecting their aggregate import onto $\mathcal{C}$ and sharing the correction equally among in-market assets: $\Pi_{\mathcal{C}}$ enforces the minimum-bid set (10) in the DA stage and additionally the increment rule (11) in the RT stage. Both are nearest-point projections of a scalar (the per-step aggregate) onto a union of intervals, available in closed form and introducing no shared integer variable.

Because the asset costs are separable and couple only through the shared limits, the game is an exact potential game with potential $\Phi = \sum_i J_i$ [14], and the proximal term makes $\Phi$ decrease monotonically each sweep. As the binaries $u_i, \delta_i$ lie in a finite set they can switch only finitely often, after which the iteration reduces to convex consensus ADMM [13]; the scheme thus converges to a mixed-integer $\varepsilon$-Nash equilibrium, with $\varepsilon$ set by the stopping tolerance on the primal and dual residuals [12].

The full scheme runs the DA GNE once per day to fix $(u_i^{\mathrm{DA}}, p_i^{\mathrm{DA}})$, locks each hour's aggregate commitment at a 65-minute gate, then runs the RT GNE every $\Delta t$, executing one step at a time.

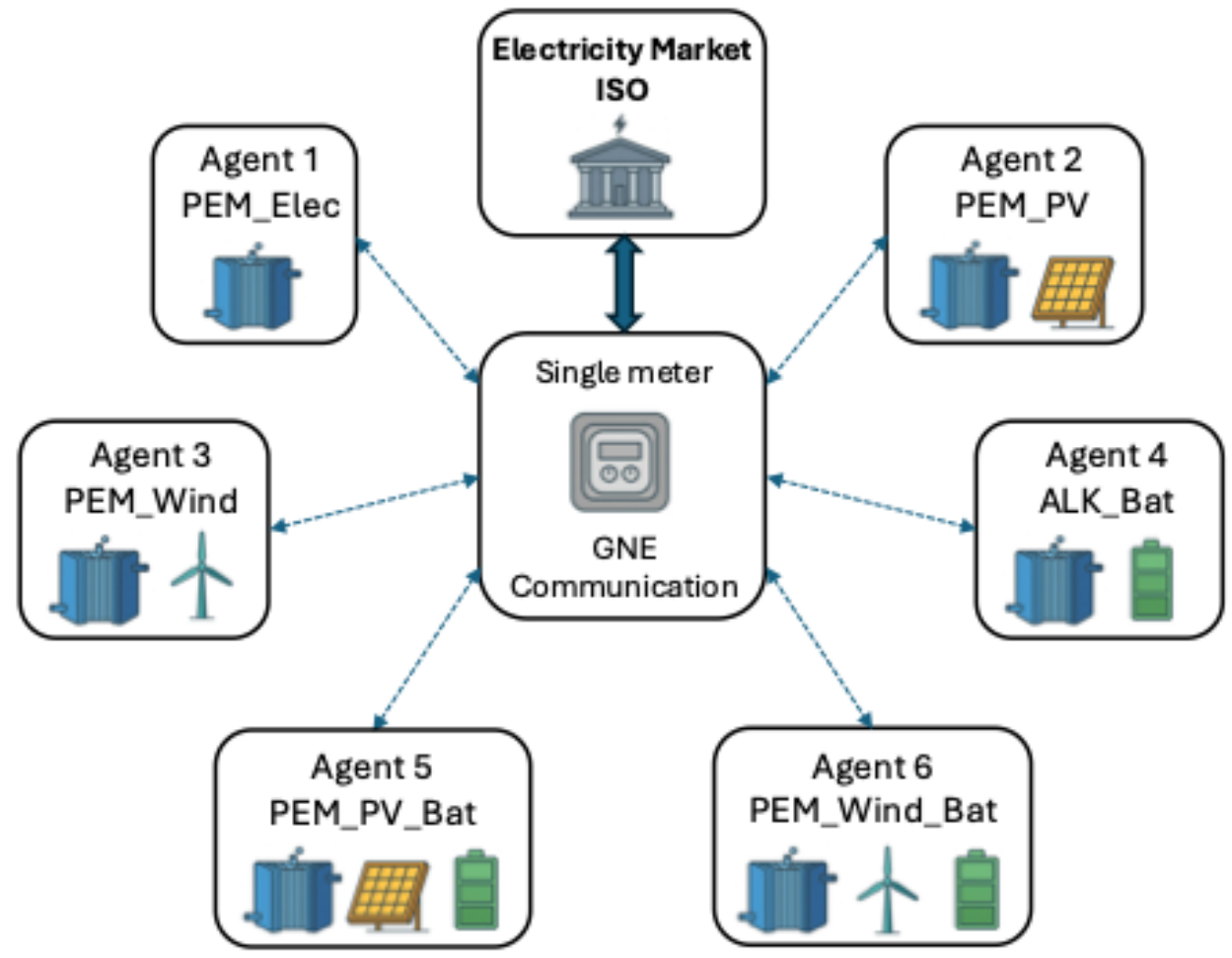

**Figure 2.** Multi-agent system architecture. Six electrochemical agents share a single grid metering point; the GNE communication channel carries only the aggregate grid import. No private cost or production data is exchanged. All agents bid jointly to the electricity market ISO and collectively produce $H_2$.

# CASE STUDY: ELECTROLYZER FLEET ON PJM

We demonstrate the framework on a heterogeneous fleet of six water electrolyzers producing hydrogen, depicted in Figure 2

and summarized in Table 1. Five agents have a proton exchange membrane (PEM) stack, and one uses an alkaline (ALK) stack; PEM units have a yield of $\eta_i = 0.020$ kg/kWh, and the ALK unit $0.018$ kg/kWh, all valued at $r_i = \$3$/kg. The units differ in rated power and in whether they carry co-located solar (PV), wind, or battery storage. Every unit has $P_i^{\max} \geq L_{\min}$, so each is individually large enough to participate alone, the setting in which the inefficiency of going it alone is most relevant.

We use real PJM Interconnection data for July 8-14, 2024. The optimizer does not see these values; it acts on forecasts obtained by random sampling around the realized series (wider in day-ahead, tighter in real time), a standard noisy-forecast model, so the receding-horizon scheme must correct forecast error online, and any forecaster can be substituted. Other settings include: $L_{\min} = 100$ kW, $L_{max} = \infty$, $\Delta t = 5$ min, RT look-ahead $H = 12$ steps (1 h), ADMM tolerance $\varepsilon = 1$ kW. All MIQPs were solved with Gurobi 13.0.1.

**Table 1:** Fleet parameters. Rated power is the maximum electrolyzer load $P_i^{\max}$; Renewable is co-located renewable capacity; Bat. is battery energy capacity. PEM = proton exchange membrane; ALK = alkaline.

| Agent | Asset | Rated power [kW] | Renewable [kW] | Bat. [kWh] |
|---|---|---|---|---|
| Agent 1 | PEM_Elec | 250 | --- | --- |
| Agent 2 | PEM_PV | 125 | 125 (PV) | --- |
| Agent 3 | PEM_Wind | 250 | 250 (Wind) | --- |
| Agent 4 | ALK_Bat | 200 | --- | 200 |
| Agent 5 | PEM_PV_Bat | 125 | 125 (PV) | 100 |
| Agent 6 | PEM_Wind_Bat | 250 | 250 (Wind) | 200 |

We compare two scenarios. In the *individual* scenario, each agent/unit participates in the power markets on its own, so the market rules apply to its own import. In the *coalition* scenario, the six agents bid jointly as a single participant, and the market rules apply only to their combined import, solved as the receding-horizon GNE above.

# RESULTS AND DISCUSSION

## Mechanism: Aggregate vs. Per-Asset Constraint

Figure 3 shows aggregate grid import across the full seven-day horizon. Under individual participation, each renewable-coupled asset must maintain $\geq 100\ kW$ of grid import individually during the daily solar peak (10:00-16:00), even when on-site generation approaches rated capacity. The GNE coalition enforces $L_{\min}$ on the aggregate only: renewable-coupled assets exit the market ($p_i = 0$), and grid-only partners cover the coalition's minimum. This selective opt-out cuts average aggregate import by 18% over the week (14% on Day, $940.8 \rightarrow 812.0$kW) and eliminates forced curtailment.

The rule also permits full market exit ($\sum_i p_i = 0$): unused this moderate-price week, but under extreme conditions (e.g. heat-wave, June 2025, peaks $> \$3700$/MWh) the coalition cuts aggregate import to zero.

## Renewable Curtailment Elimination

Figure 4 shows that fleet-wide curtailment drops from 9818 kWh (Individual) to 814 kWh (GNE), a 92% reduction over seven days. Figure 5 shows the mechanism per asset: under individual participation, the forced $\geq L_{\min}$ grid import (gray) crowds out on-site generation, which is then curtailed (red hatch); under the coalition these assets opt out, so renewable supply meets the load. PV-coupled assets drive the improvement in curtailment (Agent 2: -96%; Agent 5: -100%), since midday solar generation aligns with peak LMPs. Wind assets also benefit (Agent 6: -100%; Agent 3: -54%), though less dramatically, as wind is less temporally concentrated.

## Economic Performance

Table 2 reports seven-day financial outcomes. The GNE coalition earns \$3970 fleet profit versus \$2900 for individual participation (+37%). Electricity procurement cost falls by 29% because renewable-coupled assets exit the market during high-LMP generation peaks, consuming zero-marginal-cost on-site power instead. The 11% process output decrease is economically rational: during curtailment-forced hours the market LMP exceeds the effective product value ($\lambda_t^{RT}/1000 > r_i\eta_i$), so buying grid electricity to sustain production destroys value. The GNE coalition identifies and avoids these hours autonomously, without any central coordinator. Batteries are also far less cycled under the coalition (see Figure 5), opting out replaces the storage otherwise needed to absorb forced imports, so market flexibility substitutes for storage.

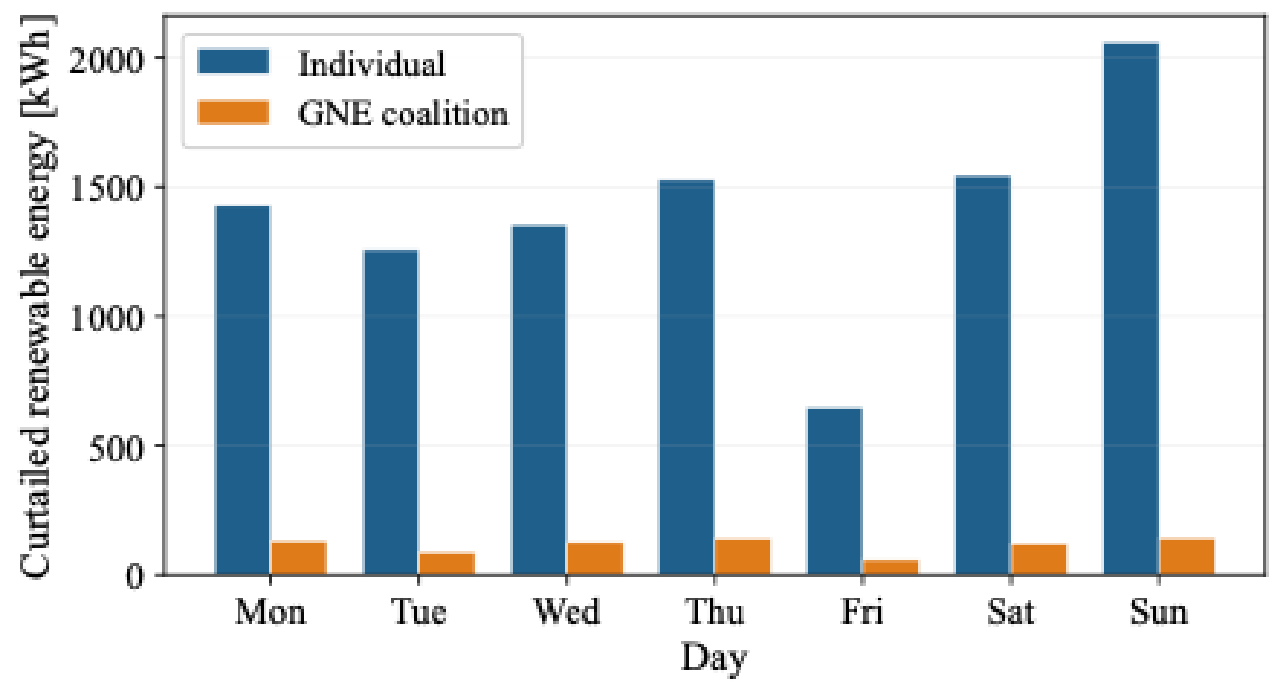

**Figure 4.** Daily fleet renewable curtailment: Individual vs. GNE coalition.

Critically, every asset earns strictly higher profit in the coalition, confirming individual rationality without side payments or central redistribution. Each asset simply solves its own local MIQP with a shared aggregate constraint; the GNE equilibrium coordinates strategies autonomously. This contrasts with aggregator-based schemes, which typically require explicit profit reallocation to ensure individual rationality [4].

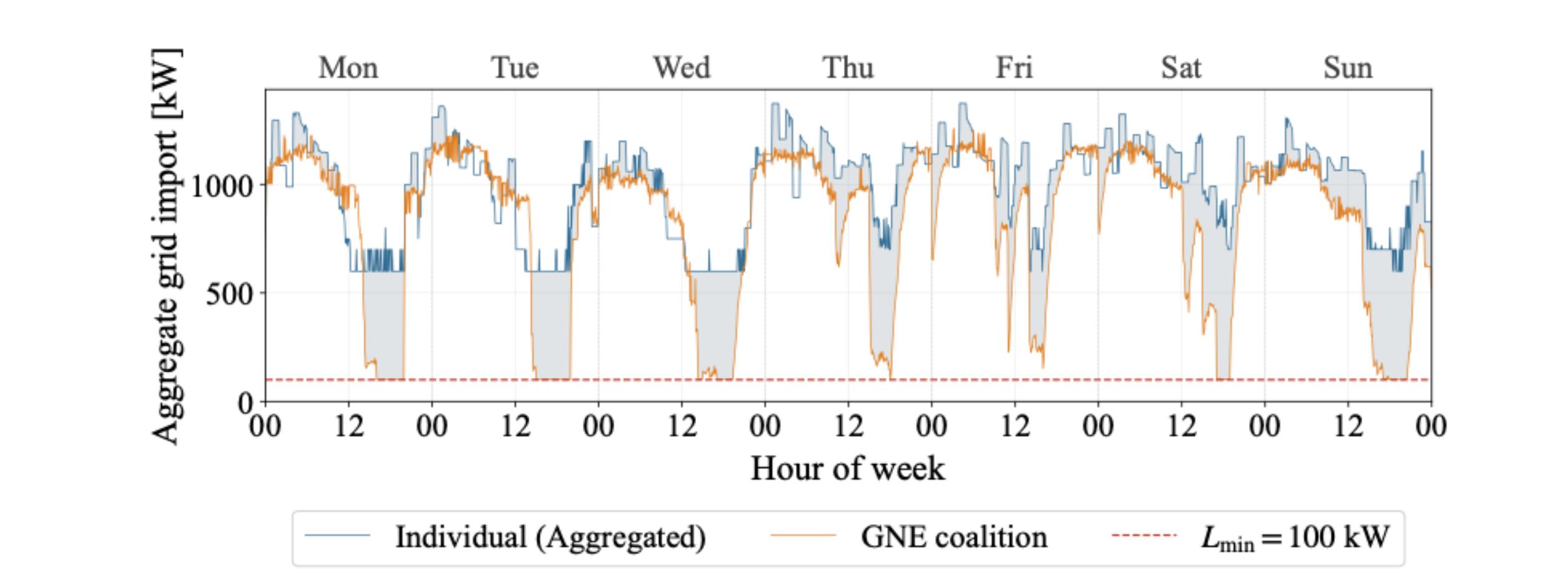

**Figure 3.** Aggregate grid import over the seven-day horizon: Individual (blue) vs. GNE coalition (orange); the shaded band is the import the coalition avoids.

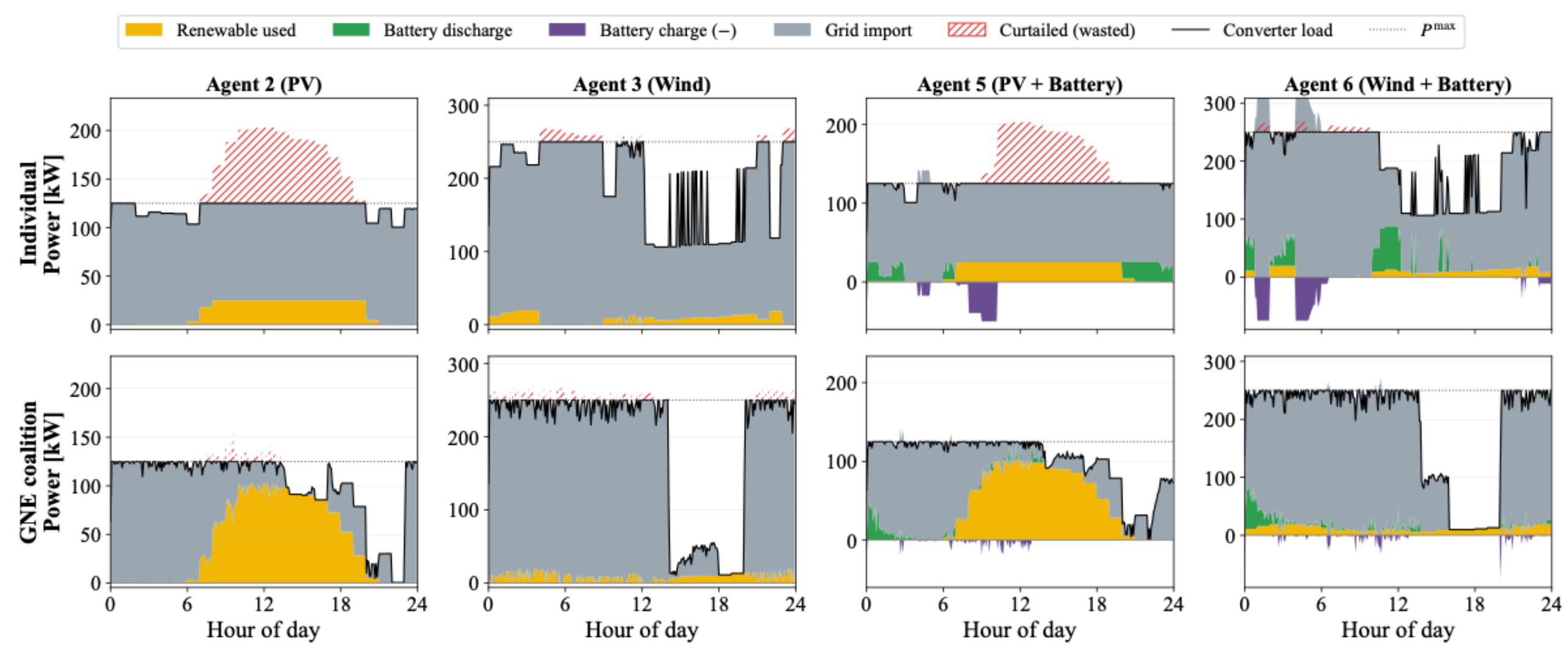

**Figure 5:** Per-agent electrolyzer power balance on a representative day (renewable-coupled agents).

**Table 2:** Seven-day simulation: Individual vs. GNE coalition.

| | Individual Participation | | | | GNE Coalition | | | |
|---|---|---|---|---|---|---|---|---|
| **Agent** | **Profit [$]** | **Cost [$]** | **$H_2$ [kg]** | **Curt. [kWh]** | **Profit [$]** | **Cost [$]** | **$H_2$ [kg]** | **Curt. [kWh]** |
| Agent 1 | 561 | 1,587 | 716 | 0 | 645 | 1,257 | 634 | 0 |
| Agent 2 | 296 | 930 | 408 | 4,258 | 665 | 409 | 358 | 167 |
| Agent 3 | 672 | 1,531 | 734 | 1,404 | 827 | 1,161 | 663 | 646 |
| Agent 4 | 266 | 1,328 | 531 | 0 | 330 | 1,034 | 455 | 0 |
| Agent 5 | 331 | 911 | 414 | 3,537 | 658 | 473 | 377 | 1 |
| Agent 6 | 774 | 1,505 | 759 | 619 | 846 | 1,188 | 678 | 0 |
| **Fleet Total** | **2,900** | **7,790** | **3,563** | **9,818** | **3,970** | **5,522** | **3,164** | **814** |
| **GNE vs. Ind.** | **+37%** | **-29%** | **-11%** | **-92%** | | | | |

## CONCLUSIONS

The receding-horizon game framework enables a self-governing coalition of electrochemical DERs to eliminate the structural curtailment inefficiency imposed by per-asset wholesale market rules. Distributed MIQP-ADMM computes the coalition GNE at each dispatch interval via closed-form aggregate projections onto the market feasibility set, with no shared binary variables, no central coordinator, and no private cost data exchanged between units. The scheme converges reliably across dispatch intervals, including under noisy day-ahead forecasts, because receding-horizon replanning absorbs forecast error online. The central novelty is embedding the two-stage day-ahead to real-time wholesale market structure within a receding-horizon GNE framework. The formulation is asset-agnostic: any electrochemical DER with a power-to-product conversion yield fits without modification. Future work will extend to stochastic renewable forecasts, larger fleets, and strategic bidding when coalition volume becomes non-negligible.

## ACKNOWLEDGEMENTS

This work is based upon work supported by the National Science Foundation under grant no. CMMI-2328160.

## DECLARATION OF USE OF AI

Claude (Anthropic) and Google's Gemini were used in the preparation of this manuscript. These tools were used to improve the grammar and fluency of the text, as well as to assist with coding and commenting.

## AUTHOR IDENTIFIERS

Brahmbhatt P: 0009-0006-1502-1359
Avraamidou S: 0000-0002-9334-9951

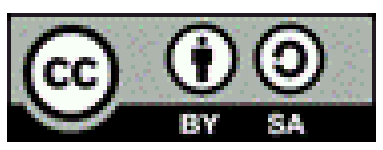